\documentclass[preprint, 12pt, english]{amsart}
\usepackage{amsthm}
\usepackage{amsmath}
\usepackage{latexsym, amssymb}
\usepackage{txfonts}
\usepackage{mathtools}
\usepackage{color}
\usepackage[all]{xy}

\newtheorem{thm}{Theorem}[section] 

\newtheorem{lem}[thm]{Lemma}

\theoremstyle{definition}
\newtheorem{rem}[thm]{Remark}

\newcommand\operA[2]{{\if!#2!\operatorname{#1}\else{\operatorname{#1}_{#2}^{\phantom{I}}}\fi}} 

\newcommand{\Trace}[1][]{\if!#1!\operatorname{Tr}\else{\operatorname{Tr}_{#1}^{\phantom{I}}}\fi} 

\long\def\forget#1\forgotten{{}} %

\def\({\left(}
\def\){\right)}

\newcommand\LAY[3][]{{\begin{array}{c}\mbox{#2} \if#1!{}\else{+}\fi \\ \mbox{#3}\end{array}}}

\makeatletter

\def\ps@pprintTitle{%
 \let\@oddhead\@empty
 \let\@evenhead\@empty
 \def\@oddfoot{}%
 \let\@evenfoot\@oddfoot}

\newcommand{\bigperp}{%
  \mathop{\mathpalette\bigp@rp\relax}%
  \displaylimits
}

\newcommand{\bigp@rp}[2]{%
  \vcenter{
    \m@th\hbox{\scalebox{\ifx#1\displaystyle2.1\else1.5\fi}{$#1\perp$}}
  }%
}
\makeatother

\renewcommand{\geq}{\geqslant}
\renewcommand{\leq}{\leqslant}

\newif\iffurther
\furtherfalse

\begin{document}

\title{Essential Dimension of Central Simple Algebras of Degree 8 and Exponent 2 in Characteristic 2}

\author{Adam Chapman}
\address{School of Computer Science, Academic College of Tel-Aviv-Yaffo, Rabenu Yeruham St., P.O.B 8401 Yaffo, 6818211, Israel}
\email{adam1chapman@yahoo.com}

\begin{abstract}
The goal of this note is to reduce the existing upper bound for the essential dimension of central simple algebras of degree 8 and exponent 2 over fields of characteristic 2 from 10 to 9.
\end{abstract}

\maketitle

\section{Introduction}

Fixing a field $k$ in the background, we consider the covariant functor $\operatorname{Alg}_{d,e}$ sending any field $F$ containing $k$ to the set of all (isomorphism classes of) central simple algebras of degree $d$ and exponent dividing $e$ over $F$. The essential dimension of an algebra $A$ in $\operatorname{Alg}_{d,e}(F)$ is the minimal possible transcendence degree of a field $E$ over $k$ with $k\subseteq E \subseteq F$ for which there exists a central simple algebra $A_0$ in $\operatorname{Alg}_{d,e}(E)$ with $A_0 \otimes F=A$.
The essential dimension of $\operatorname{Alg}_{d,e}$, denoted $\operatorname{ed}(\operatorname{Alg}_{d,e})$, is thus the supremum on the essential dimensions of $A$ where $A$ ranges over all algebras in $\operatorname{Alg}_{d,e}(F)$ and $F$ ranges over all fields containing $k$.
The question of computing $\operatorname{ed}(\operatorname{Alg}_{d,e})$ is a very difficult one (see \cite[Section 10]{Merkurjev:2013}) and the exact value is known only in specific cases: it is known that
$\operatorname{ed}(\operatorname{Alg}_{2,2})=\operatorname{ed}(\operatorname{Alg}_{3,3})=2$, but for any prime $p>3$ it is still not known what $\operatorname{ed}(\operatorname{Alg}_{p,p})$ is. The value of $\operatorname{ed}(\operatorname{Alg}_{4,2})$ is 4 when $\operatorname{char}(k)\neq 2$ and $3$ when $\operatorname{char}(k)=2$. The value of $\operatorname{ed}(\operatorname{Alg}_{8,2})$ is 8 when $\operatorname{char}(k) \neq 2$ (\cite[Corollary 1.4]{BaekMerkurjev:2012}). The exact value of $\operatorname{ed}(\operatorname{Alg}_{8,2})$  when $\operatorname{char}(k)=2$ is not known, but was bounded from above by 10 and below by 3 in \cite{Baek:2011}. The lower bound was improved to 4 in \cite{McKinnie:2017}.
The goal of this short note is to prove the following theorem:
\begin{thm}[{Main Theorem}]
When $\operatorname{char}(k)=2$, we have $\operatorname{ed}(\operatorname{Alg}_{8,2})\leq 9$.
\end{thm}
This is proved in Section \ref{Main}, and Section \ref{First} describes the main ingredients.

\section{Preliminaries}\label{First}

When $F$ is a field of $\operatorname{char}(F)=2$, a quaternion algebra over $F$ is of the form
$[\alpha,\beta)_F=F \langle i,j : i^2+i=\alpha, j^2=\beta, j i j^{-1}=i+1 \rangle$ for some $\alpha \in F$ and $\beta\in F^\times$. We write $\wp(t)$ for $t^2+t$ and $F[\wp^{-1}(\alpha)]$ thus stands for $F[i : i^2+i=\alpha]$.
By \cite[Theorem 9.1.4]{GilleSzamuely:2017} (credited to Teichm\"uller), every class in ${_2Br}(F)$ is represented by a tensor product of quaternion algebras.
In particular, when $A$ is a central simple algebra of degree 4 and exponent 2, it decomposes as a tensor product of two quaternion algebras (\cite{Racine:1974}).

A tool that makes the study of the essential dimension of central simple algebras of degree 8 and exponent 2 easier is the chain lemma for quaternion algebras (\cite[Chapter 14, Theorem 7]{Draxl:1983}):
If $[\alpha,\beta)_F=[\alpha',\beta')_F$, then there exists $\delta \in F$ for which $[\alpha,\beta)_F=[\delta,\beta)_F=[\delta,\beta')_F=[\alpha',\beta')_F$.
It is also good to have the necessary and sufficient condition for $[\alpha,\beta)_F$ to be split: $\alpha=x^2+x+y^2 \beta$ for some $x,y\in F$ (\cite[Page 104]{Draxl:1983}). In \cite{Baek:2011}, Baek was using another necessary and sufficient condition for splitting: $\beta=x^2+xy+y^2 \alpha$  for some $x,y \in F$ (\cite[Corollary 4.7.5]{GilleSzamuely:2017}), but here we prefer to use the former condition.

Another important tool is the fact that the essential dimension of compositums of $n$ cyclic quadratic field extensions of fields $F$ containing $k$ of $\operatorname{char}(k)=2$ is 1 (see \cite[Lemma 2]{Ledet:2004} and \cite[Remark 3.8]{BerhuyFavi:2003} for reference, and \cite[Remark 4.6]{ChapmanFlorenceMcKinnie:2022} for an explicit construction; despite the condition $|k|\geq 2^n$ appearing in \cite[Lemma 2]{Ledet:2004}, the proof of the latter goes through under the weaker assumption that $|F|\geq 2^n$, and for our purpose this always applies, for in order for noncommutative division algebras to exist over $F$, the cardinality of $F$ must be infinite).
The meaning is that if $K=F[\wp^{-1}(\alpha),\wp^{-1}(\beta),\wp^{-1}(\gamma)]$ is a compositum of three cyclic quadratic extensions of a field $F$ containing $k$, then there exists a field $L\supseteq k$ transcendence degree 1 over $k$ containing some elements $a,b,c$ such that the tensor product $L[\wp^{-1}(a),\wp^{-1}(b),\wp^{-1}(c)] \otimes_L F$ is $K$. In other words, the elements $\alpha,\beta,\gamma$ can be chosen to come from a single field $L$ of transcendence degree at most 1 over $k$.

The last tool is the following lemma:
\begin{lem}[{\cite[Lemma 3.3]{Baek:2011}}]\label{BaekLem}
Given a field $F$ of characteristic 2 and a field extension $E=F[i : i^2+i=\alpha]$ and $\beta \in F$, if $x^2+xy+y^2\alpha=u^2+uv+v^2\beta$ for some $x,y,u,v\in F$, then $[\beta,x+yi)_E=[\beta,y+v)_E$.
\end{lem}

\section{Proof of the main theorem}\label{Main}

Consider a division algebra $A$ of degree 8 and exponent 2 over a field $F$ containing a subfield $k$ of $\operatorname{char}(k)=2$.
By \cite{Rowen:1984}, $A$ contains a maximal subfield $K=F[\wp^{-1}(\alpha),\wp^{-1}(\beta),\wp^{-1}(\gamma)]$.
By \cite[Lemma 2]{Ledet:2004}, $\alpha$ and $\beta$ and $\gamma$ can be chosen to come from a single subfield $k \subseteq L \subseteq F$ of transcendence degree 1 over $k$.
Let $E=F[\wp^{-1}(\alpha)]=F[\mu : \mu^2+\mu=\alpha]$. Then $A_E$ is Brauer equivalent to $[\beta,b)_E\otimes [\gamma,c)_E$ for some $b,c \in E$.
The corestriction back to $F$ is trivial, and so $[\beta,N_b)_F=[\gamma,N_c)_F$, where $N_b=b_0^2+b_0b_1+b_1^2 \alpha$, and $N_c=c_0^2+c_0c_1+c_1^2\alpha$, given that $b=b_0+b_1 \mu$ and $c=c_0+c_1 \mu$ for some  $b_0,b_1,c_0,c_1 \in F$ (\cite{MammoneMerkurjev:1991}).
By the chain lemma for quaternion algebras, there exists $\delta \in F$ for which 
$[\beta,N_b)_F=[\delta,N_b)_F=[\delta,N_c)_F=[\gamma,N_c)_F$.
Therefore $\delta=\beta+u^2+u+x^2 N_b=\gamma+y^2+y+z^2 N_c=\lambda^2+\lambda+t^2 N_b N_c$ for some $u,x,y,z,\lambda,t \in F$. By replacing $\delta,y,\lambda$ with $\delta+u^2+u$, $y+u$, $\lambda+u$ respectively, we can assume $u=0$.
Note that $y=\wp^{-1}(\beta+\gamma+x^2N_b+z^2N_c)$ and $\lambda=\wp^{-1}(\beta+x^2N_b+t^2N_bN_c)$.
If $x=0$ then $\delta=\beta$, which means $[\beta+\delta,b)_E$ is split. If $x\neq 0$, then from the equation $\delta=\beta+x^2N_b$ we obtain $(\beta+\delta)\frac{1}{x^2}=b_0^2+b_0b_1+b_1^2\alpha$. It follows then from Lemma \ref{BaekLem} that $[\beta+\delta,b)_E=[\beta+\delta,b_1+\frac{1}{x})_E$.
Similarly, if $z=0$ then $[\gamma+\delta,c)_E$ is split, and otherwise $[\gamma+\delta,c)_E=[\gamma+\delta,c_1+\frac{1}{z})_E$, and if $t=0$ then $[\delta,bc)_E$ is split and otherwise $[\delta,bc)_E=[\delta,\frac{1}{t}+b_0c_1+b_1c_0+b_1c_1)_E$.

Now, $A_E=[\beta,b)_E\otimes [\gamma,c)\sim_{Br}[\beta+\delta,b)_E \otimes [\gamma+\delta,c)_E\otimes [\delta,bc)_E$. In the expression $[\beta+\delta,b)_E \otimes [\gamma+\delta,c)_E\otimes [\delta,bc)_E$, each term is either split or isomorphic to $[\beta+\delta,b_1+\frac{1}{x})_E, [\gamma+\delta,c_1+\frac{1}{z})_E$ and $[\delta,\frac{1}{t}+b_0c_1+b_1c_0+b_1c_1)_E$ respectively. In either case, each term descends to $F$.
If one of the terms is split, then $A_E$ is the restriction of a biquaternion algebra $Q_1 \otimes Q_2$ from $F$ to $E$, and thus $A=Q_1 \otimes Q_2 \otimes [\alpha,d)_F$ for some $d \in F^\times$. This means $A$ is a decomposable algebra of degree 8 and exponent 2, and by \cite[Corollary 5.9]{McKinnie:2017}, $A$ descends to a field of transcendence degree at most 4 over $k$.
We may therefore continue with the case of $x,z,t\neq 0$.
Hence, $A \sim_{Br} [\alpha,a)_F \otimes [\beta+\delta,b_1+\frac{1}{x})_F \otimes [\gamma+\delta,c_1+\frac{1}{z})_F \otimes [\delta,\frac{1}{t}+b_0c_1+b_1c_0+b_1c_1)_F$.
Consider now the algebra $B=[\alpha,a)_T \otimes [\beta+\delta,b_1+\frac{1}{x})_T \otimes [\gamma+\delta,c_1+\frac{1}{z})_T \otimes [\delta,\frac{1}{t}+b_0c_1+b_1c_0+b_1c_1)_T$ where $T=L(a,b_0,b_1,c_0,c_1,x,z,t,y,\lambda)$. Note that $y$ and $\lambda$ are algebraic over $L(a,b_0,b_1,c_0,c_1,x,z,t)$, and thus the transcendence degree of $T$ over $k$ is at most 9.
In order to conclude the argument, we need to explain why $B$ is of index 8 rather than 16, which means that it is $M_2(A_0)$ for some division algebra $A_0$ of degree 8 over $T$, and thus $A$ descends to a degree 8 division algebra over a field of transcendence degree at most 9 over $k$.
The restriction of $B$ to $R=T[\wp^{-1}(\alpha)]$ is Brauer equivalent to $[\beta+\delta,b_1+\frac{1}{x})_R \otimes [\gamma+\delta,c_1+\frac{1}{z})_R \otimes [\delta,\frac{1}{t}+b_0c_1+b_1c_0+b_1c_1)_R$.
Since $[\beta+\delta,b_1+\frac{1}{x})_R=[\beta+\delta,b)_R$ (for the same equality $\delta=\beta+x^2N_b$ holds true as before), $[\gamma+\delta,c_1+\frac{1}{z})_R=[\gamma+\delta,c)_R$ and $[\delta,\frac{1}{t}+b_0c_1+b_1c_0+b_1c_1)_R=[\delta,bc)_R$, we get that $B \otimes R$ is Brauer equivalent to $[\beta,b)_R \otimes [\gamma,c)_R$. Thus, $B$ is split by $T[\wp^{-1}(\alpha),\wp^{-1}(\beta),\wp^{-1}(\gamma)]$, and so its index is 8. \qed

\begin{rem}
It seems that in the argument in \cite{Baek:2011}, the author is suggesting that one could assume $x=z=1$ in our computation above.
If we were talking only about the quadratic equations $\delta=\beta+x^2N_b=\gamma+y^2+y+z^2 N_c$ this would be correct, because replacing $b_i$ and $c_j$ with $\frac{b_i}{x}$ and $\frac{c_j}{z}$ would eliminate their occurrences in these equations. However, that would change $[\beta,b)_E$ to $[\beta,\frac{b}{x})_E$ and $[\gamma,c)_E$ to $[\gamma,\frac{c}{z})_E$. It is therefore not clear why one can assume $x=z=1$. If this were true, that would reduce the upper bound of the essential dimension from 9 to 7.
\end{rem}

%
%
%

\section*{Acknowledgements}

The author wishes to thank the anonymous referee for the detailed reading of the manuscript and for the helpful comments and suggestions.

\bibliographystyle{alpha}
\def\cprime{$'$}

\end{document}